\documentclass[12pt]{article} 	 \usepackage{amsfonts}\usepackage{amssymb}\usepackage{amsmath}\usepackage[english]{babel}

\begin{document}
\baselineskip=14pt
\pagestyle{plain}
{\Large
 \makeatletter
\renewcommand{\@makefnmark}{}
\makeatother

\centerline {\bf On a differential operator with absent spectrum}

\medskip
\medskip
\centerline {\bf Alexander Makin}

\begin{abstract}
In this paper, we consider spectral problem for the nth  order ordinary differential operator with degenerate boundary conditions. We construct a nontrivial example of boundary value problem which has no eigenvalues.
\end{abstract}

\medskip
\medskip

\medskip
\medskip

Consider the boundary value problem for the $n$th order ordinary differential equation
$$
l(u)=u^{(n)}(x)+\sum_{m=1}^{n}p_m(x)u^{(n-m)}(x)=\lambda u(x)\eqno(1)
$$
with  boundary conditions
$$
B_j(u)=u^{(n-j)}(0)+d(-1)^{j+1}u^{(n-j)}(1)=0\eqno(2)
$$
$(j=1,\ldots n)$, where $n=2\nu$,  coefficient $d$ is an arbitrary complex number, and the complex-valued coefficients $p_m(x)$ are functions in the class
 $L_1(0,1)$. Suppose that $p_m(x)=(-1)^mp_m(1-x)$  almost everywhere on the segment $[0,1]$,
$m=1,\ldots, n$. We will study the spectrum of problem (1), (2).

 Let a function $u_k(x)$ be the solution of the Cauchy problem
$$
l(u)=\lambda u,\quad u_k^{(j)}(1/2)=\delta_{k,j},\eqno(3)
$$
where $k=0,\ldots,n-1$, $j=0,\ldots,n-1$.
Denote
$$
u_-(x)=\sum_{k=0}^{\nu-1}c_{2k+1}u_{2k+1}(x),\quad
u_+(x)=\sum_{k=0}^{\nu-2}c_{2k}u_{2k}(x),
$$
where $c_i$ are arbitrary constants $(i=0,\ldots,n-1)$.
Then
$$
u_-^{(2k)}(1/2)=0,\quad u_+^{(2k+1)}(1/2)=0,\quad k=0,\ldots,\nu-1.
$$

Obviously, that the functions $w_-(x)=-u_-(1-x)$  and $w_+(x)=u_+(1-x)$ are the solutions of equation (1) and satisfy the same initial conditions at the point $1/2$ just as the functions   $u_-(x)$ and $u_+(x)$, correspondingly. This, together with the uniqueness of the solution of Cauchy  problem (3) implies, that
$
u_-(x)=-u_-(1-x),\quad
u_+(x)=u_+(1-x)
$
for $0\le x\le1$.
It follows from this that
$$
u_-^{(n-j)}(0)+(-1)^{j+1}u_-^{(n-j)}(1)=0,\quad
u_+^{(n-j)}(0)+(-1)^{j}u_+^{(n-j)}(1)=0\eqno(4)
$$
$(j=1,\ldots, n)$. It follows from (4) that for any complex number $\lambda$ the function $u_-(x)$ is a solution of problem (1), (2) if $d=1$, and for any  complex number $\lambda$ the function $u_+(x) $ is a solution of problem (1), (2) if $d=-1$. Thus, we establish that if $d=\pm1$ the spectrum of problem (1), (2) fills all complex plane. If $p_m(x)\in C^{n-m}(0,1)$, $m=1,\ldots, n$, this assertion was proved in [1].

Assume, for a number $\lambda$ a function $\tilde u(x)$ is a solution of problem (1), (2) if $d\ne\pm1$. Then $\tilde u(x)=u_+(x)+u_-(x)$. We see that
$$
u_-^{(n-j)}(0)+u_+^{(n-j)}(0)+
(-1)^{j+1}d(u_-^{(n-j)}(1)+u_+^{(n-j)}(1))=0\eqno(5)
$$
$(j=1,\ldots, n)$.
It follows from (4), (5) that
$$
u_-^{(n-j)}(0)(1-d)+u_+^{(n-j)}(0)(1+d)=0
$$
$(j=1,\ldots, n)$. From this and the definition of the functions $u_+(x)$ and $u_-(x)$, we have
$$
(1+d)\sum_{k=0}^{\nu-2}c_{2k}u_{2k}^{(n-j)}(0)+(1-d)\sum_{k=0}^{\nu-1}c_{2k+1}u_{2k+1}^{(n-j)}(0)=0
$$
$(j=1,\ldots, n)$, hence, the constants $c_i$ $(i=0,\ldots, n-1)$ satisfy the system of linear equations
$$
\sum_{i=0}^{n-1}c_i(1+d(-1)^i)u_{i}^{(n-j)}(0)=0\eqno(6)
$$
$(j=1,\ldots, n)$. The determinant of linear system (6) is
$$
\Delta=(1-d^2)^\nu \det||u_{i}^{(n-j)}(0)||.
$$
Since the last determinant  is the Wronskian of the fundamental system of the solutions of equation (1), it is nonzero. Therefore, system (6) has only trivial solution, i.e. the function $\tilde u(x)\equiv0$. Hence, if $d\ne\pm1$ problem (1), (2) has no eigenvalues.

For the first time problem (1), (2) for $n=2$ was investigated in [2].

\medskip
\medskip
\centerline{References}
\medskip
\medskip
1. V.A. Sadovnichy, B.E. Kanguzhin. On a connection between the spectrum of a differential operator with symmetric coefficients and boundary conditions. Dokl. Akad. Nauk SSSR. 1982. V. 267. No. 2. P. 310-313.

2. M.H. Stone. Irregular differential systems of order two and the related expansion problems.  Trans. Amer. Math. Soc. 1927. V. 29. No. 1. P. 23-53.

\medskip
\medskip
\medskip
\medskip

E-mail: alexmakin@yandex.ru

}
\end{document}